\documentclass[11pt]{article}
\usepackage{graphicx}
\usepackage{epsfig}
\usepackage{latexsym}
\newsavebox{\toy}
\savebox{\toy}{\framebox[0.65em]{\rule{0cm}{1ex}}}
\newcommand{\QED}{\usebox{\toy}\end{demo}}
\newenvironment{property}%
{\begin{list}{}{\setlength{\rightmargin}{0pt}%
\setlength{\itemsep}{0pt}}}{\end{list}}
\newlength{\templength}
\newcommand{\bp}{\setlength{\templength}{\labelwidth}%
\setlength{\labelwidth}{2em}\begin{property}}
\newcommand{\ep}{\end{property}\setlength{\labelwidth}{\templength}}
\newtheorem{theorem}{Theorem}[subsection]
\newtheorem{lemma}[theorem]{Lemma}
\newtheorem{proposition}[theorem]{Proposition}
\newtheorem{corollary}[theorem]{Corollary}
\newtheorem{assumption}{Assumption}
\newtheorem{definition}{Definition}[subsection]
\newtheorem{remark}{Remark}[subsection]
\newtheorem{exercise}{Exercise}[subsection]

\newcommand{\Thm}[1]{Theorem \ref{Thm.#1}}
\newcommand{\Lem}[1]{Lemma \ref{Lem.#1}}

\newcommand{\Prop}[1]{Proposition \ref{Prop.#1}}
\newcommand{\Cor}[1]{Corollary \ref{Cor.#1}}

\newcommand{\Theorem}[1]{\begin{theorem}\label{Thm.#1}}
\newcommand{\Lemma}[1]{\begin{lemma}\label{Lem.#1}}
\newcommand{\Proposition}[1]{\begin{proposition}\label{Prop.#1}}
\newcommand{\Corollary}[1]{\begin{corollary}\label{Cor.#1}}
\newcommand{\Assumption}[1]{\begin{assumption}\label{Ass.#1}\rm}
\newcommand{\Definition}[1]{\begin{definition}\label{Def.#1}\rm}
\newcommand{\Remark}[1]{\begin{remark}\label{Rem.#1}\rm }
\newcommand{\Exercise}[1]{\begin{exercise}\label{Exe.#1}\rm }

\newcommand{\bd}{\begin{displaymath}}
\newcommand{\ed}{\end{displaymath}}
\newcommand{\bdn}{\begin{equation}}
\newcommand{\bdnl}{\begin{equation}\label}
\newcommand{\edn}{\end{equation}}
\newcommand{\barray}{\begin{array}}
\newcommand{\earray}{\end{array}}
\newcommand{\bds}{\begin{description}}
\newcommand{\eds}{\end{description}}
\newcommand{\bitemize}{\begin{itemize}}
\newcommand{\eitemize}{\end{itemize}}
\newcommand{\benumerate}{\begin{enumerate}}
\newcommand{\eenumerate}{\end{enumerate}}
\newcommand{\btabbing}{\begin{tabbing}}
\newcommand{\etabbing}{\end{tabbing}}
\newcommand{\bcenter}{\begin{center}}
\newcommand{\ecenter}{\end{center}}
\newcommand{\bflushright}{\begin{flushright}}
\newcommand{\bflushleft}{\begin{flushleft}}
\newcommand{\eflushright}{\end{flushright}}
\newcommand{\eflushleft}{\end{flushleft}}
\newcommand{\bdnn }{\begin{eqnarray*}}
\newcommand{\ednn }{\end{eqnarray*}}
\newcommand{\bdmn}{\begin{eqnarray}}
\newcommand{\edmn}{\end{eqnarray}}

\newcommand{\SSC}[1]{\section{#1}\setcounter{equation}{0}}

\newcounter{biblio}
\newenvironment{references}%
{\begin{list}{[\arabic{biblio}]}{\usecounter{biblio}%
\setlength{\leftmargin}{2.5em}\setlength{\rightmargin}{0pt}%
\setlength{\labelwidth}{2em}\setlength{\itemsep}{0pt}}}{\end{list}}
\newcommand{\References}%
{\vspace{2.8ex plus .3ex minus .3ex}%
\begin{center}{\bf References}\end{center}\begin{references}}

\usepackage{amssymb}
\usepackage{mathrsfs}
\newcommand{\bL}{{\mathbb{L}}}
\newcommand{\N}{{\mathbb{N}}}
\newcommand{\Z}{{\mathbb{Z}}}
\newcommand{\zd}{\Z^d}

\newcommand{\R}{{\mathbb{R}}}
\newcommand{\rd}{\R^d}

\newcommand{\ra }{\rightarrow }
\newcommand{\lra }{\longrightarrow }

\newcommand{\Ra}{\Rightarrow }

\newcommand{\ov}{\overline}
\newcommand{\tl}{\widetilde}

\newcommand{\Llra}{\Longleftrightarrow }

\newcommand{\vvs}{\vspace{2ex}}
\newcommand{\vs}{\vspace{1ex}}

\newcommand{\lef}{\left}
\newcommand{\rig}{\right}
\newcommand{\ri}{\right}
\newcommand{\st}{\stackrel}

\newcommand{\8}{\infty}

\newcommand{\dps}{\displaystyle}
\newcommand{\sub}{\subset}

\newcommand{\suplim}{\mathop{\overline{\lim}}}
\newcommand{\epty}{\emptyset}
\renewcommand{\a}{\alpha}
\renewcommand{\b}{\beta}
\newcommand{\gm}{\gamma}
\newcommand{\Gm}{\Gamma}

\newcommand{\del}{\delta}

\newcommand{\e}{\varepsilon}

\newcommand{\h}{\eta}

\newcommand{\kp}{\kappa}

\newcommand{\lm}{\lambda}

\newcommand{\m}{\mu}
\newcommand{\n}{\nu}
\newcommand{\rh}{\rho}

\newcommand{\s}{\sigma}

\newcommand{\vp}{\varphi}

\newcommand{\W}{\Omega}

\newcommand{\cF }{{\cal F}}

\newcommand{\cP }{{\cal P}}

\newcommand{\cR }{{\cal R}}

\newcommand{\ovZ}{\ov{Z}}
\newcommand{\ovn}{\ov{N}}
\newcommand{\ovN}{\ov{N}}

\makeatletter
\def\section{\@startsection{section}{1}{\z@}{-3.5ex plus -1ex minus 
 -.2ex}{2.3ex plus .2ex}{\bf}}
\def\subsection{\@startsection{subsection}{2}{\z@}{-3.25ex plus -1ex minus 
 -.2ex}{1.5ex plus .2ex}{\bf}}
\makeatother

\begin{document}
\parindent=0pt

\bcenter

\large{\bf Central Limit Theorem for Branching Random Walks }\\
\large{\bf in Random Environment}

\vvs \normalsize

\noindent Nobuo YOSHIDA\footnote{Partially
supported by JSPS Grant-in-Aid for Scientific
Research, Kiban (C) 17540112}\\

\ecenter

\begin{abstract}
We consider branching random walks in $d$-dimensional integer 
lattice with time-space i.i.d. offspring distributions. 
When $d \ge 3$ and the fluctuation of the environment is well moderated 
by the random walk, we prove a central limit theorem 
for the density of the population, together 
with upper bounds for the density of the most populated site and 
the replica overlap. We also discuss the phase transition of 
this model in connection with directed polymers 
in random environment. 
\end{abstract}

\vvs
\small 
Abbreviated Title: CLT for Branchig RW in Random Environment 

%%%%%%%%%%%%%%%%%%%%%%
%\small
%\tableofcontents
%%%%%%%%%%%%
\normalsize

%%%%%
\SSC{Introduction}
%%%%%%
\hspace{3mm} We consider particles in $\zd$, performing random walks 
and branching into independent copies at each step of the random walk.
When a particle occupies a site $x \in \zd$ at time 
$t \in \N=\{0,1,..\}$, then, 
it moves to a randomly choosen  adjacent site $y$ 
at time $t+1$ and is replaced by 
$k$ new particles with probability $q_{t,x}(k)$ ($k  \in \N$). 
We assume that the offspring distributions  
$q_{t,x}=(q_{t,x}(k))_{k \in \N}$ are i.i.d. in time $t$ and 
space $x$. This model was investigated earlier in \cite{Bir03,BGK05}, 
and we call it the branching random walks in random environment (BRWRE). 
See section \ref{sec:brw} for a more precise definition.

\hspace{3mm} An object of central interest in this model 
is the population $N_{t,x}$ 
of the particles at time-space $(t,x) \in \N \times \zd$,
 and the total population 
$N_t=\sum_{x \in \zd}N_{t,x}$ at time $t$. Due to the random 
environment, the population has much more fluctuation as compared with the 
non-random environment case, e.g.,\cite[section 4.2]{Rev94}.
This fluctuation results from ``disastrous locations" in 
time-space, where the offspring distribution $q_{t,x}(k)$ happens to 
assign extremely high probability to small $k$'s. 
Thanks to the random walk, on the other hand, some of the particles are 
lucky enough to elude those disastrous locations. 
Therefore, the spatial motion component 
of the model has the effect to moderate the fluctuation.

\hspace{3mm} As is discussed above, the random environment intensifies the 
fluctuation of the population, while the spatial motion moderates it.
As was observed earlier in \cite[Theorem 4]{BGK05}, these competing 
factors in the model give rise to a phase transition as follows. 
When the randomness of the offspring distribution is well moderated by 
that of the random walk, the growth of the total 
population is of the same order as 
its expectation with strictly positive probability. 
When, on the other hand, 
the randomness of the environment dominates,  
the total population grows strictly slower than its expectation 
almost surely. We also discuss this phase transition later in this 
article. Interestingly, this phase transition shares the same aspect with, 
or even explains, the localization/delocalization transition of directed 
polymers in random environment \cite{CSY04}, and of parabolic 
Anderson model with time-space i.i.d potentials, e.g.,\cite{CaMo94}. 

\hspace{3mm} In this article, we mainly 
consider the case in which the fluctuation 
caused by the random environment is well moderated by 
the random walk. It is known that this is the case 
if  $d \ge 3$ and the mean offspring is controlled by 
a sqare moment condition \cite[Theorem 4]{BGK05}--see \Thm{BCLT} below. 
We prove a central limit theorem 
for the density of the population (\Thm{BCLT}, \Cor{BCLT}), together 
with upper bounds for the density of the most populated site and 
the replica overlap (\Prop{overlap}). Our method here is based on 
sqare moment estimates. In section \ref{DP}, we discuss 
the phase transition of BRWRE in connection with that of 
the directed polymers in random environment \cite{CSY04}. 

%%%%%%%%%%%%%%%%
\subsection{Branching random walks in random environment (BRWRE)}
 \label{sec:brw}
%%%%%%

\hspace{3mm} We start with some remarks on the usage 
of the notation in this paper.
We write $\N=\{0,1,2,...\}$, 
$\N^*=\{1,2,...\}$ and 
$\Z=\{ \pm x \; ; \; x \in \N \}$. 
Let 
$(\W, \cF, P)$ be 
a probability space, which is not necessarily the one we define by
 (\ref{WcF})--(\ref{P^q}) later on.
We write $P[X]=\int X \; dP$ and 
$P[X:A]=\int_A X \; dP$ for a r.v.(random variable) $X$ and an 
event $A$. 

\hspace{3mm} We now define the model. Let $p(\cdot, \; \cdot)$ 
be a transition probability for a Markov chain 
with a countable state space $\Gm$. 
To each $(t,x) \in \N \times \Gm$, we associate a distribution 
$$
q_{t,x}=(q_{t,x}(k))_{k \in \N} \in [0,1]^\N, 
\; \; \; \sum_{k \in \N}q_{t,x}(k)=1
$$
on $\N$. Then, the branching random walk (BRW) with 
offspring distribution $q=(q_{t,x})_{(t,x) \in \N \times \Gm}$ 
is described as the following dynamics:

\bitemize
\item
At time $t=0$, there is one particle at 
the origin $x=0$. 
\item
Suppose that there are $N_{t,x}$ particles at each site $x \in \Gm$ 
at time $t$.
At time $t+1$, the $\n$-th particle at a site $x$ 
($\n =1,..,N_{t,x}$) jumps to a site 
$y=X^\n_{t,x}$ with probability $p(x,y)$ independently of each other.
 At arrival, it dies, leaving $K^\n_{t,x}$ new particles there.
\eitemize

We formulate the above description more precisely. 
%We will do so 
%in a little more general framework.  
%Our main interest is of course 
%the case when $\Gm=\zd$ and $p(\cdot, \; \cdot)$ is given by (\ref{1/2d}).
The following formulation is an analogue of \cite[section 4.2]{Rev94}, 
where non-random offspring distributions are considered. See also 
\cite[section 5]{BGK05} for the random offspring case.

\vvs
\noindent $\bullet$ {\it Spatial motion:}
A particle at time-space location $(t,x)$ is supposed to jump 
to some other location $(t+1,y)$ and is replaced by its 
children there. 
Therefore, the spactial motion should be described by assignning 
destination of the each particle 
at each time-space location $(t,x)$. So, we are guided to 
the following definition.
We define the measurable space 
$(\W_X,\cF_X)$ as the set $\Gm^{\N \times \Gm \times \N^*}$ 
with the product $\s$-field, and $\W_X \ni X \mapsto X^\n_{t,x}$ for each 
$(t,x,\n) \in \Gm \times \N \times \N^*$ as the projection. 
We define $P_X \in \cP (\W_X,\cF_X)$ 
as the product measure such that
\bdnl{P_X}
P_X (X^\n_{t,x}=y)=p(x,y)\; \; \; 
\mbox{for all $(t,x,\n) \in \N \times \Gm \times \N^*$ and $y \in \Gm$.}
\edn
Here, we interpret 
$X^\n_{t,x}$ as the position at time $t+1$ 
of the children born from the $\n$-th 
particle at time-space location $(t,x)$. 

\vvs
\noindent $\bullet$ {\it Offspring distribution:}
We set $\W_q =\cP (\N)^{\N \times \Gm}$, where 
$\cP (\N)$ denotes the set of probability measures on $\N$:
$$
\cP (\N)=\{ q=(q (k))_{k \in \N} \in [0,1]^\N
\; ; \; \sum_{k \in \N}q (k)=1\}.
$$
Thus, each $q \in \W_q$ is a function 
$(t,x) \mapsto q_{t,x}=(q_{t,x}(k))_{k \in \N}$ from 
$\N \times \Gm$ to $\cP (\N)$. We interpret $q_{t,x}$ as the 
offspring distribution for each particle which occupies the 
time-space location $(t,x)$. 
The set $\cP(\N)$ is equipped 
with the natural Borel $\s$-field induced from that of $[0,1]^\N$. 
We denote by $\cF_q$ the product $\s$-field on $\W_q$. 

We define the measurable space $(\W_K,\cF_K)$ as the set 
$\N^{\N \times \Gm \times \N^*}$ with the product $\s$-field, 
and $\W_K \ni K \mapsto K^\n_{t,x}$ for each 
$(t,x,\n) \in \N \times \Gm \times \N^*$ as the projection.
For each fixed $q \in \W_q$, we define $P^q_K \in \cP (\W_K,\cF_K)$ 
as the product measure such that
\bdnl{P^q_K}
P^q_K (K^\n_{t,x}=k)=q_{t,x}(k)\; \; \; 
\mbox{for all $(x,t,\n) \in \Gm \times \N \times \N^*$ and $k \in \N$.}
\edn
We interpret $K^\n_{t,x}$ as the number of the children 
born from the  $\n$-th particle at time-space location $(t,x)$.

We now define the branching random walk in random environment. 
We fix a product measure $Q \in \cP (\W_q,\cF_q)$, which describes the 
i.i.d. offspring distribution assigned to each time-space location. 
Finally, we define $(\W, \cF)$ by
\bdnl{WcF}
\W=\W_X \times \W_K \times \W_q, \; \; \; 
\cF=\cF_X \otimes \cF_K \otimes \cF_q, 
\edn
and $P^q, P \in \cP (\W, \cF)$ by
\bdnl{P^q}
P^q=P_X \otimes P^q_K \otimes \del_q, \; \; \; 
P=\int Q (dq)P^q.
\edn
We denote by $N_{t,x}$ the population 
at time-space location $(t,x) \in \N \times \Gm$, 
which is defined inductively
 by $N_{0,x}=\del_{0,x}$ for 
$t=0$, and 
\bdnl{n(x,t)}
N_{t,x}=\sum_{y \in \Gm}\sum_{\n =1}^{N_{t-1,y}}
\del_x (X^\n_{t-1,y})K^\n_{t-1,y}
\edn
for $t \ge 1$. The total population at time $t$ is then given by 
\bdnl{n_t}
N_t =\sum_{x \in \Gm}N_{t,x}=\sum_{y \in \Gm}
\sum_{\n =1}^{N_{t-1,y}}K^\n_{t-1,y}.
\edn
We remark that the total population is exactly the classical Galton-Watson 
process if $q_{t,x} \equiv q$, where $q \in \cP
(\N)$ is non-random.On the other hand, if $\Gm$ 
is a singleton, then $N_t$ is the polulation of the Smith-Wilkinson 
model \cite{SW69}. 

\hspace{3mm} For $p >0$, we write 
\bdmn
m^{(p)}&=& Q[m^{(p)}_{t,x}] \; \; \mbox{with} \; \; 
m^{(p)}_{t,x}=\sum_{k \in \N}k^pq_{t,x}(k), \label{m_p}\\
m &=& m^{(1)}. \label{m_1}
\edmn
Note that for $p \ge 1$, 
$$
m^p \le Q[m_{t,x}^p] \le m^{(p)} 
$$
by H\"older's inequality. 
%\bdnl{m_1}
%m=Q[m_{t,x}] \; \; \mbox{with} \; \; 
%m_{t,x}=\sum_{k \in \N}kq_{t,x}(k).
%\edn
%%%%%%%
We set 
\bdnl{ovn_t}
\ovn_{t,x}=N_{t,x}/m^t\; \; \mbox{and}\; \; \ovn_t=N_t/m^t.
\edn
$\ovn_t=N_t/m^t$ 
is a martingale (\Lem{n_t/mart} below), 
and therefore the following limit always exists:
\bdnl{ovn_8}
\ovn_\8=\lim_{t \ra \8} \ovn_t, \; \; \mbox{$P$-a.s.}
\edn
%%%%%%%%%%%%
\subsection{Results} 
\label{secBWCLT}
%%%%%%%%%%%%%%%%%%%%%%%%%%%%%%%%%%%%%%%%%%%%%%%%%%%%%%%%
Before we state our results, we fix our notation for simple random walk.

\vvs
%%%%%%%%%%%%
\noindent $\bullet$ {\it The random walk:} 
$(\{ S_t\}_{t \in \N}, P_S^x)$ is a simple random walk on 
the $d$-dimensional integer lattice $\zd$ starting from $x \in \zd$.
More precisely, we let $(\W_S, \cF_S)$ be the path space 
 $(\zd)^\N$ with the cylindrical 
$\sigma$-field, and let 
$\W_S \ni S \mapsto S_t$, $t \in \N$ be the projection. 
We define $p :\zd \times \zd \mapsto \{ 0,\frac{1}{2d}\}$ by
\bdnl{1/2d}
p(x,y)=\left\{ 
\barray{ll}
\frac{1}{2d} & \mbox{if $|x-y|=1$,} \\
0 & \mbox{if $|x-y|\neq 1$,}\earray \rig.
\edn
where $|x|=(|x_1|^2+..+|x_d|^2)^{1/2}$ for $x \in \zd$.
We consider 
the unique
probability measure $P_S^x$   on $(\W_S, \cF_S)$ such that
$S_t-S_{t-1}$, $t=1,2,..$ are independent and 
$$
P_S^x\{ S_0=x\}=1, \; \; \; 
P_S^x\{ S_t\!-\!S_{t-1}=y\}=p(0,y),\; \; \; 
\mbox{for $y \in \zd$.}
$$
In the sequel, $P_S^0$ will be simply written by $P_S$.
We define the return probability of the simple random walk:
\bdnl{pi_d}
\pi_d=P_S (S_t =0 \; \; \mbox{for some $t \ge 1$}).
\edn
As is well-known, $\pi_1=\pi_2=1$, and $\pi_d <1$ for $d \ge 3$.

To state our results, 
we assume that $\Gm=\zd$ and that 
$p(\cdot, \; \cdot)$ is given by (\ref{1/2d}). 
Then, 
with the notation introduced by (\ref{n(x,t)})--(\ref{ovn_8}), we state:
%%%%%%%%%%
\Theorem{BCLT}
%%%%%%%%
Suppose that
\bdnl{Qm_2}
m>1,\; \;m^{(2)}<\8,\; \; \mbox{and}\; \; d \ge 3.
\edn
Then, the following are equivalent:
\bds
\item[(a)] ${\dps \frac{Q[m_{t,x}^2]}{m^2} <\frac{1}{\pi_d}}$, where 
$\pi_d \in (0,1)$ is defined by (\ref{pi_d}).
\item[(b)] 
${\dps \lim_{t \ra \8} \ovN_t=\ovN_{\8}}$ in $\bL^2 (P)$.
\item[(c)]
${\dps 
\lim_{t \ra \8} \sum_{x \in \zd}\ovN_{t,x}f \lef( t^{-1/2}x\rig)
=\ovN_{\8}\int_{\rd}fg_1  }$ 
in $\bL^2 (P)$ for all $f \in C_b (\rd)$. Here and in what follows, 
\bdnl{rh_t}
g_t (x)=\lef(\frac{d}{2\pi t}\ri)^{d/2}e^{-\frac{d|x|^2}{2t}},
\; \; \; t>0,
\edn
$\int_{\rd}fg_1$ is the abbreviation for 
$\int_{\rd}f(x)g_1(x)\; dx$, 
and $C_b (\rd)$ denotes the set of bounded continuous functions on $\rd$.
\eds
%%%%%%%%
\end{theorem}
%%%%%%%
\Thm{BCLT}(a) controls the randomness of the environment in terms 
of that of the random walk. \Thm{BCLT}(b) in particular implies that 
$P (\ovN_\8 >0)>0$, i.e., the growth of the total 
population is of the same order as its expectation with 
strictly positive probability. In contrast with this, we will see that 
the total population grows strictly slower than its expectation 
almost surely, if either $d=1,2$, or the environment is random enough
(\Cor{S/E} below). It is easy to deduce from 
\Thm{BCLT}(c) the following: 
%%%%%%%%%%
\Corollary{BCLT}
%%%%%%%%
Suppose that
$$
m>1,\; \;m^{(2)}<\8,\; \; d \ge 3, 
\; \; \mbox{and}\; \; \frac{Q[m_{t,x}^2]}{m^2} <\frac{1}{\pi_d}.
$$
Then, $P (\ovN_\8 >0)>0$ and 
$$
\lim_{t \ra \8} P \lef( \lef.
\lef| \frac{1}{N_t}\sum_{x \in \zd}N_{t,x}f \lef( t^{-1/2}x\rig)
-\int_{\rd}fg_1 \rig| \ge \e  \rig| \ovN_\8 >0 \rig)=0
$$
for all $\e>0$ and $f \in C_b (\rd)$.
%%%%%
\end{corollary}
%%%%%
\Cor{BCLT} tells us that, as $t \nearrow \8$, 
the density or the spatial  distribution 
$$
\rh_{t,x}=\frac{N_{t,x}}{N_t}, \; \; \; x \in \zd 
$$
of the population converges to the standard normal distribution, 
if it is properly scaled. 
Other interesting objects related to the density would be 
$$
\rh^*_t=\max_{x \in \zd}\rh_{t,x}, \; \; \mbox{and}\; \; 
\cR_t=\sum_{x \in \zd}\rh_{t,x}^2.
$$
$\rh_t^*$ is the density at the most populated site, while 
$\cR_t$ is the probability that a given pair of particles at time $t$ are
at the same site. 
$\cR_t$ can be thought of as the replica overlap, in analogy with 
the spin glass theory. Clearly, $(\rh_t^*)^2 \le \cR_t \le \rh_t^*$. 
We use the method in this paper to show the following upper bound for $\cR_t$:
%%%%%%%%%%%
\Proposition{overlap}
%%%%%%%%%
Suppose that
$$
m>1,\; \;m^{(2)}<\8,\; \; d \ge 3, 
\; \; \mbox{and}\; \; \frac{Q[m_{t,x}^2]}{m^2} <\frac{1}{\pi_d}.
$$
Then, $P (\ovN_\8 >0)>0$ and 
$$
\cR_T= O (T^{-d/2})\; \; \; 
\mbox{in $P \lef( \cdot | \ovN_\8 >0 \rig)$-probability,}
$$
i.e., the laws 
$P \lef(T^{d/2}\cR_T \in \cdot | \ovN_\8 >0 \rig)$, $T \ge 1$ are tight.
%%%%%
\end{proposition}
%%%%%
\noindent {\bf Remarks:} After the first version of this article 
was submitted, a couple of related results are obtained. 
\bds
\item[(1)] Y. Hu and N. Yoshida \cite{HuYo07} 
prove the following localization result, 
which is in contrast with \Prop{overlap} above: 
Suppose that 
$ m^{(3)}<\8$, $Q(m_{t,x} = m) \neq 1$,  $Q(q_{t,x}(0)=0)=1$ 
and $P(\ov{N}_\8=0) =1$. 
Then, there
exists a non-random number $c \in (0,1)$ such that 
$$
\suplim_{t \nearrow \8}\cR_t
\ge c, \; \; \mbox{$P$-a.s.} 
$$ 
\item[(2)] 
Y.Shiozawa \cite{Shio07} 
considers branching Brownian motion in random environment, 
which can be thought of as a natural continuous 
counterpart of the discrete model considered in this article.
He proves \Thm{BCLT}--\Prop{overlap} for the continuous setting. 
\eds
%%%%%%%%%
\subsection{Some basic properties of $N_{t,x}$}
%%%%%%%%%%%%%%%%%%
Here again, we only assume that $(S_t,P^x_S)$ is a Markov chain 
on a countable state space $\Gm$ and with the transition probability
 $p(\cdot, \; \cdot)$. We denote the $t$ step transition probability by
\bdnl{p_t(x,y)}
p_t(x,y)=P^x_S(S_t=y).
\edn
Define $\cF_0=\{ \epty, \W\}$ and  
\bdnl{cF_t}
\cF_t =\s (X^\cdot_{s, \cdot}, K^\cdot_{s, \cdot}, 
q_{s, \cdot}\; ; \; s \le t-1 )\; \; \; t \ge 1.
\edn
This definition is natural, because the configuration of the particles 
up to time $t$ is determined by the above $\cF_t$. 
Note that 
$X^\cdot_{s, \cdot}, K^\cdot_{s, \cdot}, q_{s, \cdot}$, 
$s \ge t$ are independent of $\cF_t$.
%%%%%%%
\Lemma{Qn_t}
%%%%%%%
For $t <T$, 
\bdnl{Qn_t}
P^q[N_{T,x}|\cF_t]=\sum_{y \in \Gm}N_{t,y}
P_S^y[\prod^{T-t-1}_{u=0}m_{t+u,S_u} : S_{T-t}=x].
\edn
In particular, 
\bdnl{Qn_t2}
P^q[N_{T,x}]=P_S^0[\prod^{T-1}_{u=0}m_{u,S_u} : S_T=x]
\; \; \; \mbox{and}\; \; \; 
P^q[N_T]=P_S^0[\prod^{T-1}_{u=0}m_{u,S_u}]
\edn
%%%%%%
\end{lemma}
%%%%%%
Proof: Let $A \in \cF_t$ be arbitrary. Then, 
$$
P^q[N_{T,x}:A]
=\sum_{x_{T-1} \in \Gm}\sum_{\n \ge 0}
P^q[\del_x(X^\n_{x_{T-1,T-1}})K^\n_{T-1,x_{T-1}}: 
N_{T-1,x_{T-1}} \ge \n, A]. 
$$
By the independence, each expectation in the above sum is equal to
\bdnn
& & P_X[\del_x(X^\n_{T-1,x_{T-1}})]P_K^q[K^\n_{T-1,x_{T-1}}]
P^q[N_{T-1,x_{T-1}} \ge \n, A] \\
&= & 
p(x_{T-1},x)m_{T-1,x_{T-1}}P^q[N_{T-1,x_{T-1}} \ge \n, A].
\ednn
Hence, 
$$
P^q[N_{T,x}:A]
=\sum_{x_{T-1} \in \Gm}
P^q[N_{T-1,x_{T-1}}:A]m_{T-1,x_{T-1}}p(x_{T-1},x).
$$
By proceeding inductively, the right hand side is equal to
\bdnn
\lefteqn{\sum_{x_t,x_{t+1},..,x_{T-1}\in \Gm}P^q[N_{t,x_t}:A]
\lef( \prod^{T-1}_{u=t}m_{u,x_u}\ri)
\lef( \prod^{T-2}_{u=t}p(x_u,x_{u+1}) \ri)p(x_{T-1},x) }\\
& & = \sum_{x_t \in \Gm}P^q[N_{t,x_t}:A]
P_S^{x_t}[\prod^{T-t-1}_{u=0}m_{t+u,S_u}: S_{T-t}=x].
\ednn
Hence we have (\ref{Qn_t}). 
\hfill $\Box$
%%%%%%%
\Lemma{n_t/mart}
%%%%%%%
$(\ovn_t, \cF_t)_{t \ge 0}$ is a martingale on $(\W, \cF, P)$.
Similarly, $(P^q[\ovn_t], \cF_{q,t})_{t \ge 0}$ 
 is a martingale on $(\W_q, \cF_{q,t}, Q)$, where $\cF_{q,t}$ 
is a $\s$-field generated by $q(\cdot,\; s)$, $s \le t-1$.
%%%%%%
\end{lemma}
%%%%%%
Proof: If $t<T$, then, 
$$
P[N_T|\cF_t] 
 =  \sum_{x \in \Gm}P[N_{T,x}|\cF_t] 
 \st{\mbox{\scriptsize (\ref{Qn_t})}}{=} 
m^{T-t}\sum_{x \in \Gm}\sum_{y \in \Gm}N_{t,y}P_S^y[ S_{T-t}=x]
 = m^{T-t}N_t.
$$
\hfill $\Box$
%%%%%%%%%%%%%%%%%%%

%%%%%%%%%%%%%%%%%%%
\SSC{Proof of the results}
%%%%%%%%%%%%%%%%%%%%%%%%%%%%%%%%%%%
%%%%%%%%%%%%%%
\subsection{Lemmas}
%%%%%%%%%%
We assume that $\Gm=\zd$ and that 
$p(\cdot, \; \cdot)$ is given by (\ref{1/2d}) from here on. 
%%%%%%%%%%%%%%%%
\Lemma{N^2}
%%%%%%%%%%%%%%%%%
\bdnn
P[N_{T,x}N_{T,\tl{x}}]
&=& m^TP_S(S_T=x)\del_{x,\tl{x}} \\
& & +cm^T\sum_{t=0}^{T-1}m^tP_{S,\tl{S}}^{x,\tl{x}}
\lef[ \a^{\sum_{u=1}^t1\{S_u=\tl{S}_u\}}:S_t=\tl{S}_t, \; S_T=0\rig],
\ednn
where $\a =\frac{Q[m_{t,x}^2]}{m^2}$ and ${\dps c=\frac{m^{(2)}}{m}-1}$.
%%%%%%%%%
\end{lemma}
%%%%%%%
Proof: We follow \cite[Lemma 20]{Bir03}.
$N_{t,x}N_{t,\tl{x}}=\sum_{y,\tl{y}}F_{y,\tl{y}}$, where 
$$
F_{y,\tl{y}}=\sum_{\n =1}^{N_{t-1,y}}\sum_{\tl{\n} =1}^{N_{t-1,\tl{y}}}
K^\n_{t-1,y}K^{\tl{\n}}_{t-1,\tl{y}}
\del_x(X^\n_{t-1,y})\del_{\tl{x}}(X^{\tl{\n}}_{t-1,\tl{y}}).
$$
{\bf (1)} We first consider the expectation of $F_{y,\tl{y}}$ with
$y \neq \tl{y}$. In this case, $K^\n_{t-1,y}$ and $K^{\tl{\n}}_{t-1,\tl{y}}$ 
are independent under $P( \cdot | \cF_{t-1})$. Therefore, we have 
$$
P[F_{y,\tl{y}}|\cF_{t-1} ]=N_{t-1,y}N_{t-1,\tl{y}}m^2 p(y,x)p(\tl{y},\tl{x}).
$$
{\bf (2)} We turn to the expectation of $F_{y,\tl{y}}$ with
$y=\tl{y}$. In this case, $\{K^\n_{t-1,y}\}_{\n=1}^{N_{t-1,y}}$ 
are independent under $P( \cdot | \tl{\cF}_{t-1})$, where 
$$
\tl{\cF}_{t-1}=\s (\cF_{t-1}, (q_{t-1,x})_{x \in \zd}). 
$$
For $y = \tl{y}$ and $x = \tl{x}$, we have 
$$
P[F_{y,y}|\tl{\cF}_{t-1} ]
=N_{t-1,y}(N_{t-1,y}-1)m_{t-1,y}^2 p(y,x)^2
+N_{t-1,y}m_{t-1,y}^{(2)}p(y,x).
$$
The first and second terms on the right-hand-side 
come respectively from off-diagonal 
and diagonal terms in $F_{y,y}$.

For $y = \tl{y}$ and $x \neq \tl{x}$, we have no diagonal 
terms in $F_{y,y}$. Therefore, 
$$
P[F_{y,y}|\tl{\cF}_{t-1} ]
=N_{t-1,y}(N_{t-1,y}-1)m_{t-1,y}^2 p(y,x)p(y,\tl{x}).
$$
We now introduce the following notation:
$$
N_{t,x,\tl{x}} = N_{t,x}N_{t,\tl{x}}-N_{t,x}\del_{x,\tl{x}}.
$$
From the considerations in (1) and (2) above, and from 
$P[N_{t,x}]=m^tp_t(0,x)$, 
we obtain
\bdnn
P[N_{t,x,\tl{x}}]
&=& \sum_{y,\tl{y} \atop y \neq \tl{y}}
m^2 P[N_{t-1,y,\tl{y}}]p(y,x)p(\tl{y},\tl{x}) 
+\a m^2\sum_yP[N_{t-1,y,y}]p(y,x)p(y,\tl{x}) \\
& & +(m^{(2)}-m)\del_{x,\tl{x}}m^{t-1}p_t(0,x) \\
&=& \sum_{y,\tl{y}}P[N_{t-1,y,\tl{y}}]a(y,\tl{y})p(y,x)p(\tl{y},\tl{x}) 
+b_t (x,\tl{x}),
\ednn
where 
$$
a(y,\tl{y})=m^2\a^{1\{ y=\tl{y}\}}, 
\; \; \mbox{and}\; \; b_t (x,\tl{x})=cm^t\del_{x,\tl{x}}p_t(0,x).
$$
By \Lem{FK} below, applied to the Markov chain $(S,\tl{S})$, we get
$$
P[N_{T,x,\tl{x}}]=
c\sum_{t=0}^{T-1}m^{T+t}P_{S,\tl{S}}^{x,\tl{x}}
\lef[ p_{T-t}(0,S_t)\a^{\sum_{u=1}^t1\{S_u=\tl{S}_u\}}
:S_t=\tl{S}_t\rig].
$$
It is now easy to see that the above identity is the same as what we want.
\hfill $\Box$
%%%%%%%%%%%%%%%%
\Lemma{FK}
%%%%%%%%%%%%%%%%%
Let $S=(S_t)_{t \in \N}$ be a Markov chain with the state space 
$\Gm$. 
Suppose that 
$\vp_t$, $a_t$, $b_t$ ($t \in \N$) are functions on $\Gm$ such that
\bdnl{Sch1}
\vp_t (x)=P_S^x[a_t(S_1)\vp_{t-1}(S_1)]+b_t (x), \; \; \; 
x \in \Gm, \; \; t \ge 1,
\edn
where $P_S^x$ denotes the law of $S$, conditioned to start from $x \in \Gm$.
Then, 
\bdnl{Sch2}
\vp_T (x)=P_S^x \lef[ \vp_0 (S_T)\prod^T_{u=1}a_{T-u+1} (S_u)
+\sum_{t=0}^{T-1}b_{T-t}(S_t)\prod^t_{u=1}a_{T-u+1} (S_u)\rig], \; \; \; 
x \in \Gm, \; \; T \ge 1.
\edn
((\ref{Sch1}) and (\ref{Sch2}) are discrete analogues of 
the parabolic Schr\"odinger equation and its Feynman-Kac 
representation.)
%%%%%%%%%
\end{lemma}
%%%%%%%
Proof: Straightforward by induction on $T$. 
\hfill $\Box$
%%%%%%%%%%%%%%%%
%%%%%%%%%%
\Lemma{Nf}
%%%%%%%%%%%%%%%%%
For functions $f,\tl{f}$ on $\zd$, 
\bdnn
\lefteqn{\sum_{x,\tl{x} \in \zd}P[N_{T,x}N_{T,\tl{x}}]f(x)\tl{f}(\tl{x})}\\
&=& m^TP_S[f(S_T)\tl{f}(S_T)]+cm^T\sum_{t=0}^{T-1}m^t
P_{S,\tl{S}}^{0,0}
\lef[ \a^{\sum_{u=1}^t1\{S_t-S_u=\tl{S}_t-\tl{S}_u\}}
f(-S_T)\tl{f}(S_t-\tl{S}_t-S_T)\rig],
\ednn
where $\a =\frac{Q[m_{t,x}^2]}{m^2}$ and ${\dps c=\frac{m^{(2)}}{m}-1}$.
%%%%%%%%%
\end{lemma}
%%%%%%%
Proof: 
For $0 \le t \le T$, we compute 
\bdnn
I_t&\st{\rm def.}{=}&\sum_{x,\tl{x} \in \zd}
P_{S,\tl{S}}^{x,\tl{x}}
\lef[ \a^{\sum_{u=1}^t1\{S_u=\tl{S}_u\}}:S_t=\tl{S}_t, 
\; S_T=0\rig]f(x)\tl{f}(\tl{x}) \\
&=&\sum_{x,\tl{x} \in \zd}P_{S,\tl{S}}^{0,0}
\lef[ \a^{\sum_{u=1}^t1\{S_u-\tl{S}_u=\tl{x}-x\}}:S_t-\tl{S}_t=\tl{x}-x, 
\; S_T=-x\rig]f(x)\tl{f}(\tl{x}) \\
&=&\sum_{x,\tl{x} \in \zd}P_{S,\tl{S}}^{0,0}
\lef[ \a^{\sum_{u=1}^t1\{S_u-\tl{S}_u=S_t-\tl{S}_t\}}: 
S_t-\tl{S}_t-S_T=\tl{x}, 
\; S_T=-x\rig]f(x)\tl{f}(\tl{x}) \\
&=& P_{S,\tl{S}}^{0,0}
\lef[ \a^{\sum_{u=1}^t1\{S_t-S_u=\tl{S}_t-\tl{S}_u\}}
f(-S_T)\tl{f}(S_t-\tl{S}_t-S_T)\rig].
\ednn
By \Lem{N^2}, the left-hand-side of the desired identity equals
$$
m^TP_S[f(S_T)\tl{f}(S_T)]
+cm^T\sum_{t=0}^{T-1}m^tI_t.
$$
\hfill $\Box$
%%%%%%%%%%%%%%%%
\Lemma{SCLT}
%%%%%%%%%%%%%%%%%
Suppose that $1 \le \a <1/\pi_d$.
Then, for $f,\tl{f} \in C_b (\rd)$, 
$$
\lim_{t \ra \8}P_{S,\tl{S}}^{0,0}
\lef[ \a^{\sum_{u=0}^{t-1}1\{S_u=\tl{S}_u\}}
f(t^{-1/2}S_t)\tl{f}(t^{-1/2}\tl{S}_t)\rig]
=P_{S,\tl{S}}^{0,0}
\lef[ \a^{\sum_{u=0}^\81\{S_u=\tl{S}_u\}}\rig]
\lef(\int_{\rd}fg_1\rig)\lef(\int_{\rd}\tl{f}g_1\rig).
$$
%%%%%%%%%
\end{lemma}
%%%%%%%
Proof: This lemma is shown in the proof of 
\cite[Theorem 4.2]{Com06}.\hfill $\Box$
%%%%%%%%%
\subsection{Proof of \Thm{BCLT}}
%%%%%%%%%%%%%%%%%%
(a) $\Llra$ (b):
It follows from \Lem{Nf} that
$$
P[\ovN_t^2] 
= m^{-T}+cm^{-T}\sum_{t=0}^{T-1}m^tP_{S,\tl{S}}^{0,0}
\lef[ \a^{\sum_{u=1}^t1\{S_u=\tl{S}_u\}}\rig]
$$
and hence that
$$
\sup_{t \ge 0}P[\ovN_t^2] =\lim_{t \ra \8}P[\ovN_t^2] =
cP_{S,\tl{S}}^{0,0}\lef[ \a^{\sum_{u=1}^\81\{S_u=\tl{S}_u\}}\rig]
=cP_{S}\lef[ \a^{\sum_{u=1}^\81\{S_{2u}=0\}}\rig].
$$
The right-hand-side is finite if and only if $\a <1/\pi_d$, since 
$\sum_{u=1}^\81\{S_{2u}=0\}$ is geometrically distributed 
with the success probability $\pi_d$.

\vvs
(a) $\Ra$ (c): By a standard approximation arguments, we may assume that 
$f \in C_{b,u}(\rd)$, where $C_{b,u}(\rd)$ 
denotes the set of bounded, uniformly continuous functions.
Since (a) implies (b), it is enough to prove that, as $T \nearrow \8$, 
$$
X_T\st{\rm def.}{=}\sum_{x \in \zd}\ovN_{T,x}f \lef( T^{-1/2}x\rig)
\lra 0\; \; \; \mbox{in $\bL^2 (P)$}
$$
for $f \in C_{b,u} (\rd)$ such that  $\int_{\rd}fg_1 =0$.
By \Lem{Nf},
$$
P[X_T^2]
 =  \sum_{x,\tl{x} \in \zd}P[\ovN_{T,x}\ovN_{T,\tl{x}}]f_T(x)f_T(\tl{x})
= m^{-T}P_S[f_T(S_T)^2] +cm^{-T}\sum_{t=0}^{T-1}m^t\gm_{t,T},
$$
where $f_T(x) =f (T^{-1/2}x)$ and 
$$
\gm_{t,T}=P_{S,\tl{S}}^{0,0}
\lef[ \a^{\sum_{u=1}^t1\{S_t-S_u=\tl{S}_t-\tl{S}_u\}}
f_T(-S_T)f_T(S_t-\tl{S}_t-S_T)\rig].
$$
Since 
$$
\lim_{T \ra \8}m^{-T}\sum_{1 \le t \le T -\ln T}m^t\gm_{t,T}=0,
$$
it is enough to show that
\bds
\item[($\ast$)] \hspace{1cm}
${\dps \lim_{T \ra \8}\sup \{ \gm_{t,T}\; ; \; T -\ln T \le t <T \}=0.}$
\eds
For $T -\ln T \le t <T$, both 
$$
|T^{-1/2}S_T-t^{-1/2}S_t |  \; \; 
\mbox{and}\; \; 
|T^{-1/2}(S_t-\tl{S}_t-S_T) +t^{-1/2}\tl{S}_t| 
$$
are bounded by a constant multiple of $T^{-1/2}\ln T$. 
Since $f$ is uniformly continuous 
and $\a <1/\pi_d$, 
we have 
\bdnn
\gm_{t,T} 
&=& P_{S,\tl{S}}^{0,0}
\lef[ \a^{\sum_{u=1}^t1\{S_t-S_u=\tl{S}_t-\tl{S}_u\}}
f_t(-S_t)f_t(-\tl{S}_t)\rig]+\e_T \\
&=& P_{S,\tl{S}}^{0,0}
\lef[ \a^{\sum_{u=0}^{t-1}1\{S_u=\tl{S}_u\}}
f_t(-S_t)f_t(-\tl{S}_t)\rig]+\e_T
\ednn
with some $\e_T \ra 0$. 
Here, on the second line, 
we have used that $(S_u)_{u=0}^t\st{\rm law}{=}(S_t-S_{t-u})_{u=0}^t$.
This, together with \Lem{SCLT}, implies 
($\ast$). 

\vs
(c) $\Ra$ (b): Obvious.
\hfill $\Box$
%%%%%%%%%
\subsection{Proof of \Prop{overlap}}
%%%%%%%%
It follows from  the assumptions and \Thm{BCLT} that $P(\ov{N}_{\8}>0)>0$. 
Note that
$$
\frac{1}{N_T^2}\sum_{x \in \zd}N_{T,x}^2
=\frac{1}{\ov{N}_T^2}\sum_{x \in \zd}\ov{N}_{T,x}^2
$$
and that $\lim_{T \ra \8}\ov{N}_T=\ov{N}_{\8}>0$, 
$P(\; \cdot \; |\ov{N}_{\8}>0)$-a.s.
Therefore, it is enough to prove that 
$$
\sum_{x \in \zd}P(\ov{N}_{T,x}^2 |\ov{N}_{\8}>0)=O (T^{-d/2}).
$$
In fact, we will show that
\bds
\item[(1)] \hspace{1cm}
  $\sum_{x \in \zd}P(\ov{N}_{T,x}^2)=O (T^{-d/2})$. 
\eds 
Since $\a <1/\pi_d$, we have 
$$
0< \inf_{\zd}\Phi \le \sup_{\zd}\Phi <\8 \; \; 
\mbox{for}\; \; \Phi (x)=P^x_S [\a^{\sum_{u=1}^\81\{S_{2u}=0 \}}].
$$
Then, it follows from \Lem{Weoss} below that
\bds
\item[(2)] \hspace{1cm}
${\dps 
\sup_{x \in \zd}P_{S}^{0}
\lef[ \a^{\sum_{u=1}^t1\{S_{2u}=0 \}}:S_{2t}=x \rig]=O (t^{-d/2}),}$ 
$t \nearrow \8$.
\eds
By \Lem{N^2}, 
\bdnn
P(\ov{N}_{T,x}^2)
&= & m^{-T}P_S(S_T=x)
+cm^{-T}\sum_{t=0}^{T-1}m^tP_{S,\tl{S}}^{x,x}
\lef[ \a^{\sum_{u=1}^t1\{S_u=\tl{S}_u\}}:S_t=\tl{S}_t, \; S_T=0\rig] \\
&= & m^{-T}P_S(S_T=x)
+cm^{-T}\sum_{t=0}^{T-1}m^tP_{S,\tl{S}}^{0,0}
\lef[ \a^{\sum_{u=1}^t1\{S_u=\tl{S}_u\}}:S_t=\tl{S}_t, \; S_T=-x\rig].
\ednn
Hence, (1) can be seen via (2) as follows:
\bdnn
\sum_{x \in \zd}P(\ov{N}_{T,x}^2) 
&=& m^{-T}+m^{-T}\sum_{t=0}^{T-1}m^t
P_{S,\tl{S}}^{0,0}
\lef[ \a^{\sum_{u=1}^t1\{S_u=\tl{S}_u\}}:S_t=\tl{S}_t\rig] \\
&=& m^{-T}+m^{-T}\sum_{t=0}^{T-1}m^t
P_{S}^{0}
\lef[ \a^{\sum_{u=1}^t1\{S_{2u}=0 \}}:S_{2t}=0 \rig]=O(T^{-d/2}).
\ednn 
\hfill $\Box$
%%%%%%%%%%%%%%%%%%
\Lemma{Weoss}
%%%%%%%
Suppose that $b :\zd \ra \R$ ($d \ge 3$) is bounded and and that
$$
0< \inf_{\zd}\Phi \le \sup_{\zd}\Phi <\8 \; \; 
\mbox{with}\; \; \Phi (x)=P^x_S \lef[ 
\exp \lef( \sum_{t=0}^\8b(S_t) \ri)\ri].
$$
Then, as $T \nearrow \8$, 
$$
\sup_{x,y \in \zd}
P^x_S \lef[ \exp \lef( \sum_{t=0}^{T-1}b(S_t) \ri): S_T =y\ri]=O (T^{-d/2}).
$$
%%%%%%%
\end{lemma}
%%%%%%%%%n
Proof: Here is a receipe to get a Nash type estimate for 
Schr\"odinger semigroup, which is similar to \cite[Lemma 3.1.3]{CY04}, 
\cite[Lemma 3.3]{Va06}. 
We define 
$$
\tl{p}(x,y)={1 \over \Phi (x)}e^{b(x)}p(x,y)\Phi (y), \; \; \; 
\tl{m}(x)=e^{-b(x)}\Phi^2 (x).
$$
Then, it is easy to see that $\tl{p}$ is a transition probability of 
an $\tl{m}$-symmetric Markov chain on $\zd$. 
Moreover, we have by the assumption that
\bds
\item[(1)]\hspace{1cm} 
${\dps 0< \inf_{\zd}\tl{m} \le \sup_{\zd}\tl{m} <\8}$ 
(``strong reversibility" \cite[page 27]{Woe00}), 
\item[(2)]\hspace{1cm}
${\dps \lef( \sum_{x \in A}\tl{m}(x)\ri)^{d-1 \over d} \le \kp
\sum_{x \in A,y \not\in A}\tl{m}(x)\tl{p}(x,y)}$, 
for any finite $ A \sub \zd$, \\
where the constant $\kp$ is independent of $A$ 
(``$d$-isoperimetric inequality" \cite[page 40]{Woe00}),
\item[(3)]\hspace{1cm}
${\dps \sum_{y \in \zd \atop |x-y| \le r}\tl{m}(y) \ge \e r^d}$ for 
any $r \in \N$ and $x \in \zd$, \\
where the constant $\e>0$ is independent of $r$ and $x$.
\eds
Then, by \cite[page 148, Corollary 14.5]{Woe00}, these imply that 
\bds
\item[(4)] \hspace{1cm} 
${\dps \sup_{x,y \in \zd}\tl{p}_{T}(x,y)=O(T^{-d/2})}$ as 
$T \nearrow \8$, \\
where $\tl{p}_{T}$ is the $T$ step transition 
function obtained from $\tl{p}$.
\eds 
Since 
$$
P^x_S \lef[ \exp \lef( \sum_{t=0}^{T-1}b(S_t) \ri): S_T =y\ri]=
\Phi (x)\tl{p}_{T}(x,y){1 \over \Phi(y)},
$$
the lemma follows from (4). \hfill $\Box$
%%%%%%%%%%%%%%%%%%
%%%%%%%
\SSC{Relation to directed polymers in random environment} \label{DP}
%%%%%%%%%%%%%%%%%%%
We now relate the BRWRE with directed polymers in random environment.
%%%%%%%%%%%%
\subsection{Directed polymers in random environment (DPRE)} 
%%%%%%%%%%%%%%%%%%%%%%%%%%%%%%%%%%%%%%%%%%%%%%%%%%%%%%%%

$\bullet$ {\it The random environment:}  
$\h =\{\h_{t,x} : (x,t) \in \zd \times \N \}$ is a 
sequence of r.v.'s which are real valued,
non-constant, and i.i.d. r.v.'s 
 defined on a probability space $(\W_\h, \cF_\h, Q)$ such that 
\begin{equation} \label{expint}
Q[\exp (\b \h_{t,x})] <\8 \; \; \; \mbox{for all $\b \in \R$.}
\end{equation}
We define 
\begin{equation} \label{lm(b)}
\lm (\b)=\ln Q[\exp (\b \h_{t,x})].
\edn
%%%%%%%%

$\bullet$ {\it The polymer measure:}
For any $T \in \N^*$, 
define the probability measure $\m_T$ on the path space 
$(\W_S, \cF_S)$ by  
\begin{equation} \label{mnen}
d\m_T =\frac{1}{Z_T}
\exp (\b H_T)  \;dP_S,
\end{equation}
where $\b >0$ is a parameter (the inverse temperature),
\begin{equation} \label{Ham}
H_T=\sum_{t=0}^{T-1} \h_{t,S_t} \; \; \; \mbox{and}\; \; \; 
Z_T=P_S\lef[\exp \lef(\b H_T\rig) \rig]
\end{equation}
are the Hamiltonian and the normalizing constant (the partition function).

Define the {\it normalized partition
  function} by 
\begin{equation} \label{Wn}
\ovZ_T=Z_T \exp (-T \lm (\b ))\;, \; \; T \geq 1. 
\end{equation}
This is a 
 positive mean-one $(\cF_{\h,T})$-martingale 
on $(\W_\h, \cF_\h, Q)$ with $\cF_{\h,0}=\{\epty, \W_\h\}$ and 
$$
\cF_{\h,T}=\s (\h (\cdot, s)\; ; \; s \le T-1), \; \; \; T \ge 1.
$$
%%%%%%%%%%%%%%
% Indeed,
%for all fixed path $S$, $H_T (S, \eta)$, $T \ge 1$
%is a random walk defined on $(\W_\h, \cF_\h, Q)$, and 
%\bdmn
%\bz_T= \exp \Big( \b H_T-T\lm (\b )\Big),
%\label{en}
%\edmn
%is its exponential martingale. Hence, the average over $S$, 
% $\ovZ_T=P_S[\bz_T]$,  is itself   a 
%mean-one, positive $(\cF_{\h,T})$-martingale on $(\W_\h, \cF_\h, Q)$.
%%%%%
%As we will see now, the martingale property makes
%this sequence much easier to study than $\ln Z_T$ itself, a fact which
%was used first by Bolthausen in \cite{Bol89}.
%%%%%%%%%%%%%%%%%
By the martingale convergence theorem, the limit 
\begin{equation} \label{W8}
\ovZ_\8=\lim_{T \ra \8} \ovZ_T
\end{equation}
exists $Q$-a.s. Moreover, there are only two possibilities for the 
positivity of the limit;
\begin{equation} \label{Z8>0}
Q\{\ovZ_\8>0 \}=1\;,
\end{equation}
or 
\begin{equation} \label{Z8=0}
Q\{\ovZ_\8=0 \}=1\;.
\end{equation}
Indeed, 
the event $\{ \ovZ_\8=0\}$ is in the tail $\s$-field:
$$
\bigcap_{t \geq 1}\s [ \h (\cdot,s) \; ; \; s \geq t ]\;.
$$ 
By Kolmogorov's zero-one law, every event in the tail 
$\s$-field has 
probability 0 or 1. 

The above contrasting situations 
(\ref{Z8>0}) and (\ref{Z8=0}) will be called  the {\bf weak disorder}
phase and the {\bf strong disorder} phase, respectively. 
%%%%%%%%%%%%
\subsection{BRWRE and its associated DPRE} 
%%%%%%%%%%%%%%%%%%%%%%%%%%%%%%%%%%%%%%%%%%%%%%%%%%%%%%%%
Suppose that we are given an environment 
$\h_{t,x}$ for the directed polymer. 
We can then associate an environment, i.e., an i.i.d. random 
offspring distribution $q_{t,x}$ for the BRWRE so that 
\bdnl{e=m}
e^{\b \h_{t,x}}=m_{t,x}=\sum_{k \in \N}kq_{t,x}(k)
\edn
holds.
Among many ways to do so, let us take:
\bdnl{q=Poi}
q_{t,x}(k)=e^{-m_{t,x}}\frac{m_{t,x}^k}{k!}.
\edn
Then, by \Lem{Qn_t},
\bdnl{DP=BRW}
P^q[N_{t,x}]=P_S[\exp (\b H_T): S_T=x]
\; \; \mbox{and}\; \; \; 
P^q[N_T]=Z_T.
\edn
These imply
$$
P^q[\ovn_{T,x}]=P_S[\exp (\b H_T-\lm (\b )T): S_T=x],
\; \;  
P^q[\ovn_T]=\ovZ_T
$$
and 
$$
\m_T (S_T=x)=\frac{P^q[N_{t,x}]}{P^q[N_T]}
=\frac{P^q[\ovN_{t,x}]}{P^q[\ovN_T]}.
$$
Since $\m_T$ and $\ovZ_T$ are invariant 
under constant addition to $\h_{t,x}$, we may assume that
$$
m=Q[m_{t,x}]=Q[\exp (\b \h_{t,x})]>1
$$
without loss of generality. Moreover, by (\ref{q=Poi}) and our integrability 
assumption (\ref{expint}) on $\h_{t,x}$, 
$$
\sum_{k \in \N}k^2q_{t,x}(k)=m_{t,x}+m_{t,x}^2 \in \bL^1 (Q).
$$
Therefore, \Thm{BCLT} implies the following central limit theorem 
for DPRE, which is 
a weaker version of the results obtained in 
\cite{Bol89,SoZh96,CY06}. 
%%%%%%%%%%%%%%%%%%%
\Corollary{DPCLT}
%%%%%%%%%%%%%%%%%%%%%
 If $d \ge 3$ and $\lm (2\b)-2\lm (\b) < \ln (1/\pi_d)$, then 
\bdnl{DPCLT}
 \lim_{T \ra \8}\m_T [f( T^{-1/2}S_T)]=\int_{\rd} fg_1 
\; \; \; \mbox{in $Q$-probability for any $f \in C_b(\rd )$,}
\edn
%%%%%%
\end{corollary}
%%%%%%%%%%%

%%%%%%%%%%%%
\subsection{Phase transitions of BRWRE and DPRE} 
\label{secS/E}
%%%%%%%%%%%%%%%%%%%%%%%%%%%%%%%%%%%%%%%%%%%%%%%%%%%%%%%%
As we mentioned before, BRWRE undergoes the following phase transition. 
When the randomness of the offspring distribution is well moderated by 
that of the random walk (as in \Cor{BCLT}, \Prop{overlap}), the 
growth of the total 
population is of the same order as 
its expectation with strictly positive probability.
When, on the other hand, 
the randomness of the environment dominates,  
the total population grows strictly slower than its expectation 
almost surely. We now relate this phase transition with that for 
DPRE. 
%%%%%%%
\Proposition{S/E}
%%%%%%%%%%%%%%%
Suppose that an environment $\h_{t,x}$ for the DPRE  
and an offspring distribution $q_{t,x}$ for the BRWRE
are related so that (\ref{e=m}) holds. Then,
\bds
\item[(a)] 
$P[\ovn_\8] \le Q[\ovZ_\8]$. In particular, 
$Q(\ovZ_\8 =0)=1$ (strong disorder for DPRE) implies 
$P(\ovn_\8 =0)=1$ (the total population grows strictly 
slower than its expectation almost surely).
\item[(b)] 
The converse to (a) is not true.
\eds
%%%%%%
\end{proposition}
%%%%%%%%%%
Proof:
(a): We have:
$$
P[\ovn_\8]=\int Q(dq)P^q[\ovn_\8] 
\le \int Q(dq)\lim_t P^q[\ovn_t]
=Q[\ovZ_\8],
$$
where the inequality follows from Fatou's lemma.\\
(b) As is mentioned in \cite[Theorem 4]{BGK05},
It follows from a comparison with the Galton-Watson model 
that $P(N_\8 =0)=1$ and hence 
$P(\ovn_\8 =0)=1$ as soon as $m \le 1$. 
We have always $m \ge 1$  when $\h_{t,x}$ is a negative r.v. Meanwhile, 
the polymer is in weak disorder phase for 
$d \ge 3$ and $\b>0$ small enough.
\hfill $\Box$

\vvs
By \Prop{S/E}, we can translate the 
results from DPRE \cite{CSY03,CoVa06} to 
the following observation for BRWRE:
%%%%%%%%%%%%%%%%%%
\Corollary{S/E}
%%%%%%%%%%%%%%%%%%
Suppose one of the following conditions: 
\bds \item[(a1)] $d=1, \;
Q(m_{t,x} = m) \neq 1.$ 
\item[(a2)] $d=2, \; Q(m_{t,x} = m) \neq 1.$ 
\item[(a3)] ${\dps d \ge3, \;
Q\lef[\frac{m_{t,x}}{m}\ln\frac{m_{t,x}}{m}\ri]>\ln (2d)}$. 
\eds
Then, $P(\ovn_\8 =0)=1$. Moreover, in cases (a1) and (a3), there
exists a non-random number $c>0$ such that 
\bdnl{N_t<m^t}
\suplim_{t}\frac{\ln \ov{N}_t }{t}< -c, \; \; \mbox{a.s.} 
\edn
%%%%%%%%
\end{corollary}
%%%%%%%%%
\Cor{S/E} says that the total population grows strictly slower
than its expectation almost surely, in low dimensions or in
``random enough" environment. The result is in contrast with the
non-random environment case, where $P(\ovn_\8 =0)=1$ only for
offspring distributions with very heavy tail, more precisely, if
and only if $P [K_{t,x}^\n \ln K_{t,x}^\n]=\8$ \cite[page 24,
Theorem 1]{AtNe72}. Here, we can have $P(\ovn_\8 =0)=1$ even when
$K_{t,x}^\n$ is bounded. Also, (\ref{N_t<m^t}) is in sharp
contrast with the non-random environment case, where it is well
known --see e.g., \cite[page 30, Theorem 3]{AtNe72} --that
$$
\{ N_\8 >0\}\st{\mbox{\scriptsize a.s.}}{=} \{ \lim_{t}\frac{\ln
\ov{N}_t }{t}=0\}\; \; \; \mbox{whenever $m >1$.}
$$
\Cor{S/E} is also 
stated in \cite[Theorem 4]{BGK05}, however without (\ref{N_t<m^t}).

\hspace{3mm} We have discussed the relation between BRWRE and DPRE. 
We remark that the branching Brownian motions in random environment 
recently considered by Shiozawa \cite{Shio07} and 
Brownian directed polymers \cite{CY05} are in similar relation. 
Therefore, by the results in \cite{CY05}, \Prop{S/E} and 
\Cor{S/E} (except (\ref{N_t<m^t}) in case (a1)) 
can be extended to branching Brownian motions in random environment. 
%%%%%%%%%%%%%%%%%%%%%%%%%%%%%%%%%%%%%
\SSC{Survival and extinction}
%%%%%%%%%%%%%%%%%%%%%%%%%%%%%%
\hspace{3mm} We now close this article with 
a brief discussion on survival/extinction 
for this model, i.e., ${\bf e}\st{\rm def}{=}P(N_\8 =0)<1$ or $=1$. 
By standard computations of 
generating function, we see that 
${\bf e}^{\rm GW} \le {\bf e} \le {\bf e}^{\rm SW}$, 
where ${\bf e}^{\rm GW}$ and ${\bf e}^{\rm SW}$ stands respectively for 
extinction probabilities 
for Galton-Watson model with offspring distribution 
$Q[q_{t,x}(\; \cdot \; )]$ and Smith-Wilkinson model \cite{AtNe72,SW69}.
On the other hand, by using oriented percolation, one can construct 
examples for ${\bf e}^{\rm GW}<{\bf e}=1$ and for 
${\bf e} < {\bf e}^{\rm SW}=1$ \cite{Fuk07}.

\vvs
\small
{\bf Acknowledgements:}
The author thanks Francis Comets for drawing his attention to \cite{Com06}. 
He also thanks Ryoki Fukushima, Yuichi Shiozawa, and Shinzo Watanabe 
for discussions. 
%%%%%%%%%%%%%%%%%%

%%%%%%%%%%

\vvs
Division of Mathematics \\
Graduate School of Science \\
Kyoto University,\\
Kyoto 606-8502, Japan.\\
email: {\tt nobuo@math.kyoto-u.ac.jp}\\
URL: {\tt http://www.math.kyoto-u.ac.jp/}$\widetilde{}$ {\tt nobuo/}
%%%%%%%%%%%%%%%%%%%%%%%
\end{document}